\documentclass [12pt] {article}

\usepackage {amssymb}
\usepackage {exscale}

\title {A Fedosov Star Product of Wick Type\\
        for K\"ahler Manifolds
        }

\author {{\bf M. Bordemann\thanks{mbor@phyq1.physik.uni-freiburg.de}~,
            \addtocounter{footnote}{2}
            S. Waldmann\thanks{waldman@phyq1.physik.uni-freiburg.de}}\\[3mm]
             Fakult\"at f\"ur Physik\\Universit\"at Freiburg\\
          Hermann-Herder-Str. 3\\79104 Freiburg i.~Br., F.~R.~G\\[3mm]
        }

\date{FR-THEP-96/9 \\[1mm]
      May 1996 \\[1mm]{\small Revised Version}\\[5mm]
      }

\newcommand {\W} {\mathcal W}
\newcommand {\WL} {\mathcal W \otimes \Lambda}

\newcommand {\degs} {{\rm deg}_s}
\newcommand {\dega} {{\rm deg}_a}
\newcommand {\degh} {{\rm deg}_\hbar}
\newcommand {\Deg} {{\rm Deg}}

\newcommand {\WEYL} {{\mbox{\rm \tiny Weyl}}}
\newcommand {\WICK} {{\mbox{\rm \tiny Wick}}}
\newcommand {\weyl} {\circ}
\newcommand {\wick} {\circ'}
\newcommand {\adwick} {{\rm ad}_\WICK}
\newcommand {\adweyl} {{\rm ad}_\WEYL}

\newcommand {\BEQ} [1] {\begin {equation} \label {#1}}
\newcommand {\EEQ} {\end {equation}}

\newcommand {\proof} {{\sc Proof: }}
\newcommand {\QED} {\hfill $\square$}
\newcommand {\cc} [1] {\overline {{#1}}}

\newtheorem {LEMMA} {Lemma} [section]
\newtheorem {PROP} [LEMMA] {Proposition}
\newtheorem {THEOREM} [LEMMA] {Theorem}

\begin {document}
\maketitle

\begin {abstract}
In this letter we compute some elementary properties of the Fedosov
star product of Weyl type, such as symmetry and order of differentiation.
Moreover, we define the notion of a star product of Wick type on every
K\"ahler manifold by a straight forward generalization of the corresponding
star product in $\mathbb C^n$: the corresponding sequence of bidifferential
operators differentiates its first argument in holomorphic directions
and its second argument in antiholomorphic directions. By a Fedosov type
procedure we give an existence proof of such star products for any K\"ahler
manifold.
\end {abstract}

\newpage

\section {Introduction}

The concept of deformation quantization has been defined by
Bayen, Flato, Fronsdal, Lichnerowicz and Sternheimer in 1978
\cite{BFFLS78}:
A formal star product on a symplectic manifold $(M,\omega)$ is defined 
as a formal local associative deformation $\ast$ of 
the commutative algebra of complex-valued smooth functions  
$C^\infty (M)$ on $M$ where the formal parameter is identified with $\hbar$. 
These data
are equivalent to giving a sequence $M_t$, $t\geq 0$, of locally 
bidifferential operators on $M$ which satisfy the following axioms for
any three $f,g,h\in C^\infty (M)$:
\begin{eqnarray}
    f*g  & = & \sum_{t=0}^\infty(i\hbar/2)^tM_t(f,g)\\
   M_0(f,g) & = & fg \\
   M_1(f,g)-M_1(g,f) & = & \{f,g\} \\
   0                 & = & \sum_{t+u=s}\big( M_t(M_u(f,g),h)-
                                              M_t(f,M_u(g,h)\big) 
                                {\rm ~~for~all~}s\geq 0 \\
   M_t(1,f)=M_t(f,1) & = & 0 {\rm ~~for~all~}t\geq 1 \\
   \overline{M_t(f,g)} & = & (-1)^tM_t(\bar{g},\bar{f}){\rm ~~for~all~}
                                                  t\geq 0
\end{eqnarray}                                
In most of the literature the axiom $M_t(g,f)=(-1)^tM_t(f,g)$ is added,
leading to what we shall call a star product of Weyl type because in
flat $\mathbb R^{2n}$ with the standard Poisson structure it exactly
corresponds to the Weyl symmetrization rule in canonical quantization.
The existence proofs of DeWilde-Lecomte \cite{DL83} and Fedosov
\cite{Fed86}
refer to this kind of star products. However, in Berezin's 1974 paper
\cite{Ber74} on quantization via coherent states on K\"ahler manifolds
another type of star products appears where the operators $M_t$ are complex
and have the property that one of their arguments is differentiated in 
holomorphic directions whereas the other arguments is differentiated in 
antiholomorphic directions only (which one does strongly depend on the sign
conventions applied). See also \cite{CGR II}, \cite{BMS94}, and 
\cite{BBEW96} for more details and examples of such star products which
we shall call star products of Wick type because in flat $\mathbb C^n$
they correspond to the Wick ordering (or normal ordering) rule.

The aim of this Letter is twofold: firstly we are going to prove
some properties of the operators $M_t$ occurring in the star products
constructed by Fedosov \cite{Fed94}: in particular we prove the
(no doubt folklore) statement that these star products are of Weyl type
and that the operators $M_t$ are bidifferential of order $t$ (Lemma
\ref{FedWeylprefinal}, Theorem \ref{FedWeylfinal}).

Secondly, motivated by a discussion with Fedosov we transfer his procedure
of constructing star products to K\"ahler manifolds by which we get
a relatively simple direct proof of existence of star products of Wick type 
on each such
manifold (Theorem \ref{FedWickfinal}). Our proof does not make use of
analytic techniques usually employed in the literature to obtain 
star products as an asymptotic expansion of, say, a geometric quantization
scheme.

\section {Some Notation} \label {NotationSec}

We will use Einstein's summation convention (i.e. summation 
over repeated indices is understood) and make use of the
notation of Fedosov as in \cite {Fed94} with the following exeptions:
We denote a connection by the symbol $\nabla$ instead of Fedosov's 
$\partial$ and use symmetric tensor fields and insertion maps $i_s$
instead of Fedosov's $y$'s and derivations with respect to $y$.

Let $(M, \omega)$ be a $2n$-dim. symplectic manifold with symplectic
form $\omega = \frac{1}{2} \omega_{ij} dx^i \wedge dx^j$
where $(x^1, \ldots, x^{2n})$ are local coordinates. The Poisson
tensor corresponding to $\omega$ is defined by
$\Lambda = \frac{1}{2} \Lambda^{ij} \partial_i \wedge \partial_j$ where
$\partial_i = \partial/\partial x^i$ and
$\Lambda^{ij} \omega_{kj} = \delta^i_k$ according
to the sign convention of \cite {AM85}. 
Then the Poisson bracket of two
functions $f,g$ is given by $\Lambda (df, dg)$.

Next one defines the Fedosov algebra $\WL$ by
\BEQ {FedAlgDef}
    \WL := \left({\mathsf X}_{s=0}^\infty \mathbb C\left(
           \Gamma\left(\mbox{$\bigvee$}^s T^*M
           \otimes \mbox{$\bigwedge$} T^*M\right)\right)\right)[[\hbar]]
\EEQ
and hence elements $a\in \WL$ are of the form
$a = \sum_{r,s=0}^\infty \hbar^r a_{rs}$ where $a_{rs} \in
\Gamma(\mathbb C(\bigvee^s T^*M \otimes \bigwedge T^*M))$ are smooth
sections. Then in $\WL$ a fibrewise product is declared by the
$\vee$-product of symmetric forms in the first factor
and the $\wedge$-product of antisymmetric forms in the second factor. For
factorizing elements $f\otimes \alpha$ and $g\otimes \beta$ this reads
$(f\otimes \alpha)(g\otimes \beta) =
f\vee g \otimes \alpha \wedge \beta$.
Clearly $\WL$ is a formally $\mathbb Z \!\times\! \mathbb Z$-graded algebra
with respect to this product in the sense that the symmetric degree map
$\degs$ and the antisymmetric degree map $\dega$ of forms are
derivations.
We denote by $i_s (X)a$ the insertion (symmetric substitution)
of a vector field $X$ in the symmetric part of $a$ and by
$i_a(X)a$ the usual insertion (inner product) in the
antisymmetric part.
Then the fibrewise Weyl product is defined by
\BEQ {FibWeylDef}
    \begin {array} {c}
    \displaystyle
    a \weyl b := \sum_{r=0}^\infty \left(\frac{i\hbar}{2}\right)^r
                \Lambda^{(r)} (a,b) \\
    \displaystyle
    \Lambda^{(r)} (a,b) := \frac{1}{r!}
    \Lambda^{i_1j_1} \cdots \Lambda^{i_rj_r}
    i_s(\partial_{i_1}) \cdots i_s (\partial_{i_r}) a \;
    i_s(\partial_{j_1}) \cdots i_s (\partial_{j_r}) b
    \end {array}
\EEQ
where $a,b \in \WL$. This product is known to be an
associative deformation of the fibrewise product.
Denote by $\degh$ the degree map with respect to the order of the
formal parameter $\hbar$ (i.e. $\degh(\hbar^k a) = k \hbar^k a$)
and define the {\em total degree} by $\Deg := 2\degh + \degs$.
Then $(\WL, \weyl)$
is a formally graded algebra with respect to $\Deg$ and $\dega$ again in
the sense that $\Deg$ and $\dega$ are derivations with respect to
the fibrewise Weyl product.
The notation $a^{(k)}$ stands for an element in the Fedosov algebra with
total degree $k$ and $a^{(k)}_s$ stands for an element with total degree $k$ and
symmetric degree $s$.
We define the graded $\weyl$-commutator by
$[a,b]_\WEYL := a\weyl b - (-1)^{kl} b \weyl a$ where $k,l$ are the
antisymmetric degrees of $a$ and $b$. Moreover we define
$\adweyl(a) := [a,\cdot ]_\WEYL$ which is a graded derivation of
antisymmetric degree $k$ if $\dega a = ka$. We denote by $\W$
the elements in $\WL$ with vanishing antisymmetric degree
and by $\WL^a$ the elements with antisymmetric degree $a$.
Denote by $C$ the complex conjugation $C a := \cc a$
and by $P_\hbar := (-1)^{\degh}$ the $\hbar$-parity. Then clearly
$C$ and $P_\hbar$ are $\dega$-graded involutions
of the algebra $\WL$ with respect to $\weyl$, i.e. for $a,b \in \WL$
with $\dega a = ka$, $\dega b = lb$ we have
$C (a \weyl b) = (-1)^{kl} b \weyl a$ and analogously for
$P_\hbar$.

Now let $M$ be a K\"ahler manifold. In a holomorphic chart
$(z^1, \ldots, z^n)$ the symplectic form can be written as
$\omega = \frac{i}{2} \omega_{k\cc l} dz^k \wedge d\cc z^l$ with a
positive definite Hermitian matrix $\omega_{k\cc l}$. The Poisson tensor
is given by $\Lambda = \frac{2}{i} \omega^{k\cc l} Z_k \wedge \cc Z_l$
where $\omega^{k\cc l} \omega_{t \cc l} = \delta^k_t$ is the inverse
matrix and $Z_k = \partial/\partial z^k$
and $\cc Z_l = \partial/\partial \cc z^l$ are local base fields of type
$(1,0)$ and $(0,1)$. Since $M$ is K\"ahler we can characterize the forms
by their type and we will denote by $\pi^{(k,l)}_s$ the projection
on symmetric forms of type $(k,l)$
(see \cite [II chap. 9, 10] {KN69} for details). Moreover,
in addition to the fibrewise Weyl product we can define
a fibrewise Wick product for $a,b \in \WL$ by
\BEQ {FibWickDef}
    \begin {array} {c}
    \displaystyle
    a \wick b := \sum_{r=0}^\infty \left(\frac{i\hbar}{2}\right)^r
                 {\Lambda'}^{(r)} (a,b) \\
    \displaystyle
    {\Lambda'}^{(r)} (a,b) :=
                 \frac{1}{r!} \left(\frac{4}{i}\right)^r
                 \omega^{k_1\cc l_1} \cdots \omega^{k_r\cc l_r}
                 i_s(Z_{k_1}) \cdots i_s (Z_{k_r}) a \;
                 i_s(\cc Z_{l_1}) \cdots i_s (\cc Z_{l_r}) b .
    \end {array}
\EEQ
Then $(\WL, \wick)$ is again formally graded with respect to
$\Deg$ and $\dega$ and we define graded $\wick$-commutators
analogously as the graded
$\weyl$-commutators and set $\adwick (a) := [a, \cdot]_\WICK$ which is
again a graded derivation with respect to the fibrewise Wick product.
On a K\"ahler manifold the fibrewise Weyl and Wick product are fibrewise
equivalent in the following sense:
\begin {LEMMA} \label {WeylWickEquiLem}
For $a,b \in \WL$ we have
\BEQ {WWEqui}
    a \weyl b \; = \; S^{-1} \left( Sa \wick Sb \right)
\EEQ
where $S := e^{\hbar\Delta}$ and 
$\Delta := \omega^{k\cc l} i_s (Z_k) i_s (\cc Z_l)$.
\end {LEMMA}

If $*$ is a star product for the symplectic manifold $M$ then
\BEQ {StarMrDef}
    f * g = \sum_{r=0}^\infty \left(\frac{i\hbar}{2}\right)^r M_r (f,g)
\EEQ
for $f,g \in C^\infty (M)$ and the $M_r$ are some bidifferential operators.
The star product is said to be of {\em Weyl type} iff
$M_r (f,g) = (-1)^r M_r (g,f)$ and the $M_r$ are real.
If $M$ is a K\"ahler manifold a differential operator
$L: f \to L(f)$ is called of type $(1,0)$ resp. of type $(0,1)$ iff
in a holomorphic chart the function $f$ is
only differentiated in holomorphic resp. antiholomorphic directions.
This characterization is clearly independent of the holomorphic chart.
Then a star product on a K\"ahler manifold is said to be
of {\em Wick type} iff $M_r$ is of type $(1,0)$ in the
first argument and of type $(0,1)$ in the second argument
for all $r \ge 1$.

\section {The Fedosov star product of Weyl type}

We shall now briefly recall Fedosov's construction of a star product for an
arbitrary symplectic manifold (see \cite{Fed94}). First we have to define 
some
maps. Let $a \in \WL$ with $\degs a = ka$, $\dega a = la$ and let
$\nabla$ be an (always existing) symplectic torsion-free connection with
curvature tensor $R$.
\BEQ {deltaDef}
    \delta a := (1 \otimes dx^i) i_s (\partial_i) a 
    \qquad \qquad
    \delta^{-1} a := \left\{ \begin {array} {l}
                   0 \; \mbox{ \rm if } k+l = 0 \\
                   \frac{1}{k+l} (dx^i \otimes 1)
                   i_a (\partial_i) a \;
                   \mbox { \rm else}
                   \end {array}
                   \right. 
\EEQ
\BEQ {nablaRDef}
    \begin {array} {c}
    \nabla a := (1 \otimes dx^i) \nabla_{\partial_i} a \\
    R := \frac{1}{4} \omega_{it} R^t _{jkl} dx^i \vee dx^j \otimes
         dx^k \wedge dx^l \in \WL
    \end {array}
\EEQ
Then $\delta$ and $\nabla$ turn out to be graded derivations
with respect to the fibrewise Weyl product and
we have $\delta^2 = (\delta^{-1})^2 = 0$,
$\nabla^2 = \frac{i}{\hbar} \adweyl (R)$,
$\nabla\delta + \delta\nabla = 0$ and
$\delta \delta^{-1} a + \delta^{-1} \delta a = a$
if $\degs a + \dega a \ne 0$. Moreover Fedosov has proved the following
two theorems:
\begin {THEOREM} \label {FedWeylTheoI}
There exists a unique section $r \in \WL^1$ with $\delta^{-1} r = 0$
and $\delta r = R + \nabla r + \frac{i}{\hbar} r \weyl r$. It can be
calculated recursively with respect to the total degree $\Deg$ by
\BEQ {WeylrRecus}
    \begin {array} {c}
    r^{(3)} = \delta^{-1} R \\
    \displaystyle
    r^{(k+3)} = \delta^{-1} \left( \nabla r^{(k+2)} + \frac{i}{\hbar}
                \sum_{l=1}^{k-1} r^{(l+2)} \weyl r^{(k-l+2)} \right) .
    \end {array}
\EEQ
Then the Fedosov derivation
$D := - \delta + \nabla + \frac{i}{\hbar} \adweyl (r)$ has square zero:
$D^2 = 0$.
\end {THEOREM}
Since $D$ is a graded derivation the kernel of $D$ is a
$\weyl$-subalgebra and one defines
\BEQ {WKern}
    \W_D := \ker D \cap \W.
\EEQ
\begin {THEOREM} \label {FedWeylTheoII}
Let $D$ be constructed as in the last theorem. Then $\W_D$ is
in bijection with $C^\infty (M)[[\hbar]]$ and the projection
$\sigma: \W_D \to C^\infty (M)[[\hbar]]$ onto the part with symmetric
degree zero is a $\mathbb C[[\hbar]]$-linear isomorphism. The inverse 
$\tau = \sigma^{-1}$ for a
function $f \in C^\infty (M)$ can be constructed recursively with respect
to the total degree $\Deg$ by
\BEQ {WeyltauRecurs}
    \begin {array} {c}
    \tau (f)^{(0)} = f \\
    \displaystyle
    \tau (f)^{(s+1)} = \delta^{-1} \left(\nabla \tau(f)^{(s)}
                       + \frac{i}{\hbar} \sum_{t=1}^{s-1}
                       \adweyl \left(r^{(t+2)}\right) \tau(f)^{(s-t)}
                       \right).
    \end {array}
\EEQ
Then a star product for $M$ is given by
$f*g := \sigma (\tau(f) \weyl \tau(g))$.
\end {THEOREM}
Since $f*g = \sigma (\tau(f) \weyl \tau(g))$ is a star product
it can be written in the form (\ref{StarMrDef}) with some
bidifferential operators $M_r$.
Although the following results may generally be believed to be true they do
nevertheless not seem to have appeared in the literature:
\begin {LEMMA} \label{FedWeylprefinal}
\begin {enumerate}
\item $C r = r$ and $P_\hbar r = r$ hence $r$ is real and
      depends only on $\hbar^2$.
\item $ C \tau = \tau C$ and $P_\hbar \tau = \tau P_\hbar$
      and hence $\tau (f)$ depends only on $\hbar^2$ for
      $f\in C^\infty (M)$.
\item $M_s (f,g) = (-1)^s M_s (g,f)$ and hence $*$ is of the Weyl type.
\item $C (M_s (f,g)) = M_s (C f, C g)$
      and hence the $M_r$ are real.
\item Complex conjugation is an antilinear $*$-involution:
      $\cc{f*g} = \cc g * \cc f$.
\end {enumerate}
\end {LEMMA}
\begin {THEOREM} \label{FedWeylfinal}
\begin {enumerate}
\item For all nonnegative integers $k,l$ with $l \le [k/4]$ and 
      $k-4l\geq 1$ if $k\geq 1$ the map
      $f\mapsto\tau(f)^{(k)}_{k-4l}$ is a differential operator of order
      $k-2l$. 
\item The Fedosov star product is a Vey star product, i.e. the
      bidifferential operator $M_s$ is of order $s$ in both arguments
      and we have the formula
      \BEQ {WeylMr}
      M_s (f,g) = \sum_{k=0}^{\left[\frac{s-1}{2}\right]}
                  \sum_{l=0}^k (-4)^k
                  \Lambda^{s-2k} \left(
                  \tau(f)^{(s-2k+4l)}_{s-2k}, \tau(g)^{(s+2k-4l)}_{s-2k}
                  \right) .
      \EEQ
\end {enumerate}
\end {THEOREM}
\proof This theorem is proved by lengthy but straight forward
induction using the preceding lemma in particular the facts that $r$
and $\tau(f)$ depend only on $\hbar^2$ for $f\in C^\infty (M)$
as well as the recursion formulas for $r$ and $\tau$.
\QED

\section {The Fedosov star product of Wick type}

Now we consider a K\"ahler manifold $M$ and investigate the fibrewise
Wick product. For the map $\nabla$ we will always use the K\"ahler
connection. Then we notice the following properties of the fibrewise
Wick product using lemma \ref {WeylWickEquiLem}:
\begin {PROP} \label {DiverseProp}
For $a,b \in \WL$ we have:
\begin {enumerate}
\item $\adwick (a) = S \circ \adweyl (S^{-1} a) \circ S^{-1}$.
\item $\adwick (a) = 0 \; \Longleftrightarrow \; \adweyl (a) = 0
       \; \Longleftrightarrow \; \degs a = 0$
\item In the first order of $\hbar$ the graded fibrewise Weyl and the
      graded fibrewise Wick commutators are the same.
\item $[\delta, \Delta] = [\delta, S]
      = \left[\delta, S^{-1}\right] = 0$.
\item $\delta = - \frac{i}{\hbar} \adwick (\tilde{\delta}) 
              = - \frac{i}{\hbar} \adweyl (\tilde{\delta})$ 
      with $\tilde{\delta} = \omega_{k\cc l} dz^k \otimes d\cc z^l + 
      \omega_{k\cc l} d\cc z^l \otimes dz^k$ and hence $\delta$ is a
      graded derivation of antisymmetric degree $1$ with respect to the
      fibrewise Wick product.
\item $[\nabla, \Delta] = [\nabla, S] = \left[\nabla, S^{-1}\right] = 0$.
      Furthermore $\nabla$ is a graded derivation of antisymmetric
      degree $1$ with respect to the fibrewise Wick product and
      we have
      \[
          \nabla^2 \; = \; \frac{i}{\hbar} \adwick (SR)
          \; = \; \frac{i}{\hbar} \adwick (R)
      \]
      with $R = \frac{1}{2} i \omega_{k\cc t} R^{\cc t} \,_{\cc l i \cc j}
      dz^k \vee d\cc z^l \otimes dz^i \wedge d\cc z^j$ being the same
      element in $\WL$ as in (\ref{nablaRDef}).
\end {enumerate}
\end {PROP}
The proofs can be done easily in a holomorphic chart using the properties 
of the K\"ahler connection. Note that $SR = R + \hbar\Delta R$ with
$\degs\Delta R = 0$ and hence $\adwick(\Delta R) = 0$.

The following two theorems can be proved in a manner completely analogous to
the original theorems of Fedosov (Theorem \ref{FedWeylTheoI}
and \ref{FedWeylTheoII}) if we use the properties of the maps
$\delta$, $\nabla$ and the preceding proposition \cite {Wal95}.
\begin {THEOREM} \label {FedWickTheoI}
There exists a unique section $r' \in \WL^1$ with
$\delta^{-1} r' = 0$ and
$\delta r' = R + \nabla r' + \frac{i}{\hbar} r' \wick r'$. It can be
calculated recursively with respect to the total degree $\Deg$ by
\BEQ {WickrRecurs}
    \begin {array} {c}
    {r'}^{(3)} = \delta^{-1} R \\
    \displaystyle
    {r'}^{(k+3)} =
    \delta^{-1} \left( \nabla {r'}^{(k+2)} +
    \frac{i}{\hbar}
    \sum_{l=1}^{k-1} {r'}^{(l+2)} \wick {r'}^{(k-l+2)} \right) .
    \end {array}
\EEQ
Then the Fedosov derivation
$D' := - \delta + \nabla + \frac{i}{\hbar} \adwick (r')$ has square zero:
${D'}^2 = 0$.
\end {THEOREM}
\begin {THEOREM} \label {FedWickTheoII}
Let $D'$ be constructed as in theorem \ref{FedWickTheoI}.
Then $\W_{D'}$ is again in bijection to $C^\infty (M)[[\hbar]]$
and the projection $\sigma: \W_{D'} \to C^\infty (M)[[\hbar]]$
is again an isomorphism. The inverse $\tau' = \sigma^{-1}$ for a
function $f \in C^\infty (M)$ can be constructed recursively
with respect to the total degree $\Deg$ by
\BEQ {WicktauRecurs}
    \begin {array} {c}
    \tau'(f)^{(0)} = f \\
    \displaystyle
    \tau'(f)^{(s+1)} = \delta^{-1} \left(\nabla \tau'(f)^{(s)}
                       + \frac{i}{\hbar} \sum_{t=1}^{s-1}
                       \adwick \left({r'}^{(t+2)}\right) \tau'(f)^{(s-t)}
                       \right).
    \end {array}
\EEQ
Then a star product for $M$ is given by
$f*'g := \sigma (\tau'(f) \wick \tau'(g))$. The maps $\tau'(\cdot)^{(k)}$
are differential operators of order $k$.
\end {THEOREM}
We shall call this star product the 
{\em Fedosov star product of Wick type} on $M$. The first two bidifferential
operators of this star product are easily seen to take the following form:
\[
    f *' g \; = \; fg + i\hbar \Lambda(\partial f, \cc \partial g)
    + \cdots .
\]
We will now prove that this star product is indeed a star product of
Wick type as defined in section \ref {NotationSec}: to this end we will
examine the section $r'$ and the corresponding maps 
$\tau'$ in some more detail:
\begin {LEMMA} \label {CCWickLem}
With the notations from above we have:
\begin {enumerate}
\item $C r' = r'$ and hence
      $C {D'a} \; = \; D' C a \;$ for $a \in \WL$.
\item $C \tau' (f) \; = \; \tau' (C f)\;$
      for $f \in C^\infty (M)[[\hbar]]$.
\item $\cc{f*'g} \; = \; \cc g *' \cc f \;$ and hence the
      complex conjugation is a $*'$-involution too.
\item For all nonnegative integers $k,l$ with $l\leq [k/2]$ and
      $k-2l\geq 1$ for all $k\geq 1$ the
      linear map $f\mapsto\tau'(f)^{(k)}_{k-2l}$ is a differential
      operator of order $k-l$.
\end {enumerate}
\end {LEMMA}
\proof
Using the recursion formulas for $r'$ and $\tau'$ this lemma is
proved by a straight forward induction.
\QED

\begin {LEMMA}
With the notations form above we have for the section $r'$ and
any $p \ge 0$
\BEQ {rGrade}
    \pi^{(0,p)}_s r' \; = \; \pi^{(p,0)}_s r' \; = \; 0 .
\EEQ
\end {LEMMA}
\proof
The crucial point for this lemma is the choice of $R$ as starting point
in the recursion formula (\ref{WickrRecurs}) and not
$SR = R + \hbar \Delta R$ since $\degs \Delta R = 0$. Then it is proved
again by induction.
\QED

\begin {LEMMA}
Let $f \in C^\infty (M)$ be a holomorphic function in an
open set $U \subseteq M$ and smooth outside of $U$.
Then we have for $0 < p  \in \mathbb N$
\BEQ {TauGradeHol}
    \left. \pi^{(0,p)}_s \tau'(f) \right| U = 0 .
\EEQ
If on the other hand $f$ is antiholomorphic in $U$ then we have
\BEQ {TauGradeAnti}
    \left. \pi^{(p,0)}_s \tau'(f) \right| U = 0 .
\EEQ
\end {LEMMA}
\proof
Again we prove this lemma by induction using the preceding lemma and
the recursion formula for the $\tau'(f)$.
\QED

Now we can prove that the star product constructed in theorem
\ref{FedWickTheoII} is indeed a star product of Wick type. Note that in
the proof the associativity of $*'$ is crucial.
\begin {THEOREM} \label{FedWickfinal}
Let $*'$ be the star product constructed as in theorem
\ref{FedWickTheoII}. Then the bidifferential operators ${M'}_r$ 
are of order $r$ and hence $*'$ is a Vey star product. Moreover 
we have
\BEQ {MrLambdaWick}
  {M'}_r (f,g) = \sum_{k=0}^s\sum_{l=0}^k(-2{\bf i})^k
                  \Lambda'^{(s-k)}\left( \tau'(f)^{(s-k+2l)}_{s-k},
                                     \tau'(g)^{(s+k-2l)}_{s-k} \right)
\EEQ
and the operators ${M'}_r$ are of type $(1,0)$ in the first
and of type $(0,1)$ in the second argument for all $r$. Hence the star
product is of Wick type. For a function $f$ holomorphic on the open set
$U \subseteq M$ and for $g$ antiholomorphic on $U$ and $h$ arbitrary
we have
\BEQ {hWickHol}
    (h *' f)|U \; = \; hf|U \qquad (g *' h)|U \; = \; gh|U .
\EEQ
\end {THEOREM}
\proof
The first part is again proved by an analogous induction as 
in theorem \ref{FedWeylfinal}. 
Secondly we prove equation (\ref{hWickHol}). 
Since this statement is local
we can work in a holomorphic chart $(z^1, \ldots, z^n)$. Let $f$ be a
holomorphic function in the domain $U$ of this chart. Then for $p>0$ we
have $\pi^{(0,p)}_s \tau'(f)|U = 0$ according to the last lemma.
This implies that
\[
    \left. \pi^{(0,0)}_s
           \omega^{k_1 \cc l_1} \cdots \omega^{k_r \cc l_r}
           i_s (Z_{k_1}) \cdots i_s(Z_{k_r}) \tau' (h) \;
           i_s (\cc Z_{l_1}) \cdots i_s(\cc Z_{l_r}) \tau' (f)
    \right|U = 0
\]
for $r > 0$. Hence in $h *' f = \pi^{(0,0)}_s (\tau'(h) \wick \tau'(f))$
only the lowest order $r=0$ of the fibrewise Wick product contributes and
this implies that $h *' f = hf$. The antiholomorphic case is proved
analogously.

In order to prove that all ${M'}_r$ are of type $(1,0)$ in
the first argument we use
induction on $r$. For ${M'}_0$ and ${M'}_1$ this is obviously true so let
us assume that ${M'}_0, \ldots, {M'}_{r-1}$ are all of type $(1,0)$
in the first argument. Since ${M'}_r$ is of order $r$ it can be written as
\[
    {M'}_r (f,g) = \sum_{{|I| + |J| \le r \atop |K| + |L| \le r}}
    M^r_{IJKL}
    \frac{\partial^{|I|+|J|} f} {\partial z^I \partial \cc z^J}
    \frac{\partial^{|K|+|L|} g} {\partial z^K \partial \cc z^L}
\]
with some smooth locally defined functions $M^r_{IJKL}$. Now let $p$ be a
point such that for some multi-indices $I_0, J_0, K_0, L_0$ we have
$M^r_{I_0J_0K_0L_0} (p) \ne 0$. Then we can adjust our
holomorphic chart such that $p=0$ in this chart. Now we consider
polynomials in the coordinates $z^{K_0} \cc z^{L_0}$  and
$\cc z^{J_0} *' z^{I_0} = \cc z^{J_0} z^{I_0}$ since $z^{I_0}$ is locally
holomorphic and we use equation (\ref{hWickHol}).
By the associativity of $*'$ we get the equation
\BEQ {Dummy}
   (\cc z^{J_0} z^{I_0}) *' (z^{K_0} \cc z^{L_0})
   = \cc z^{J_0} *' (z^{I_0} *' (z^{K_0} \cc z^{L_0}))
   = \cc z^{J_0} (z^{I_0} *' (z^{K_0} \cc z^{L_0}))
\EEQ
since $\cc z^{J_0}$ is locally
antiholomorphic. Then at $z=0$ the term proportional to $\hbar^r$ 
of $(\cc z^{J_0} z^{I_0}) *' (z^{K_0} \cc z^{L_0})$ is equal to 
$cM^r_{I_0J_0K_0L_0} (0) \ne 0$ with some 
positive combinatorical factor $c$.
But this leads to a contradiction if $|J_0| >0$ since then 
the right hand side in (\ref{Dummy}) vanishes at $0$ but 
the term proportional to $\hbar^r$ of the left hand 
side is equal to $cM^r_{I_0J_0K_0L_0}(0) \ne 0$
at $z=0$. Hence only holomorphic derivatives can occur in the first argument 
of ${M'}_r$. The statement about its second argument is proved analogously.
\QED

\section* {Acknowledgements}

The authors would like to thank B. V. Fedosov for a discussion in which
he suggested that the Wick case can be treated by his methods.

\end {document}